\documentclass[12pt]{article}
\usepackage{amsmath,epsfig,amssymb,float,graphicx,amsfonts, amsbsy, xspace,colortbl}

\begin{document}

\title{{\bf Covariance fields}}

\author{\Large Nikolay Balov
\footnote{Florida State University, Department of Statistics}\\
{\small\texttt{balov@stat.fsu.edu}}}

\maketitle

\begin{abstract}
We introduce and study covariance fields of distributions on a Riemannian manifold.
At each point on the manifold, covariance is defined to be a symmetric and positive definite (2,0)-tensor.
Its product with the metric tensor specifies a linear operator on the respected tangent space. 
Collectively, these operators form a covariance operator field. 
We show that, in most circumstances, covariance fields are continuous. 
We also solve the inverse problem: 
recovering distribution from a covariance field. Surprisingly, 
this is not possible on Euclidean spaces. 
On non-Euclidean manifolds however, covariance fields are true distribution representations.
\end{abstract}

{\section{Preliminaries}}

Subject of this study are random variables on Riemannian manifolds. For the sake of clarity and 
self-consistency we will briefly recall the main notations and facts from Riemannian geometry we are going to use.
For a comprehensive introduction the reader is suggested to refer to \cite{carmo-riemannian},  
\cite{lee-smooth} or \cite{chavel-riemannian}.

{\subsection{Riemannian manifolds}}

Let M be a n-manifold with differentiable structure given as a collection of charts 
$(U_\alpha,{\bf x}_\alpha)$ where $U_\alpha$ are open sets in $\mathbb{R}^n$ and ${\bf x}_\alpha : U_\alpha \to M$ are injective.
For $p \in {\bf x}_\alpha(U_\alpha)$, $(U_\alpha,{\bf x}_\alpha)$ is called a parametrization or 
system of coordinates at p. 
Thus, when we say coordinates $x$ at a point of M, 
we will understand a local system of coordinates given by a chart $(U,{\bf x})$.

With $M_p$ we denote the tangent space at $p\in M$. The tangent bundle $TM$ on M is given by 
$TM=\{(p,v)| p\in M, v\in M_p\}$. It is a 2n-manifold. 
The map $\pi:TM \to M$, $\pi(p,v)=p$ denotes the natural projection.

Recall that if $M_1$ and $M_2$ are two manifolds and $\phi: M_1 \to M_2$ is a differentiable mapping, 
differential of $\phi$ at $p\in M_1$ (also called a push-forward) 
is a linear mapping $d\phi_p :(M_1)_p \to (M_2)_{\phi(p)}$ given by 
$d\phi_p(v) (f) = v(f\circ \phi)$ for any $v\in (M_1)_p$ and $f\in C^{\infty}(M_1)$, 
a differentiable function on $M_2$. If $d\phi_p$ is an isomorphism, 
then $\phi$ is a local diffeomorphism at p (Theorem 2.10 in \cite{carmo-riemannian}).

With respect to a parametrization $(U,{\bf x})$ and $p={\bf x}(x_1,...,x_n)\in {\bf x}(U)$, 
the tangent space $M_p$ of M at p has canonical basis $\{\frac{\partial}{\partial x_i}|_p\}_{i=1}^n$.
Let $v\in M_p$, then $v_x=(v_x^1,...,v_x^n)\in\mathbb{R}^n$, such that $v=\sum_{i=1}^n v_x^i \frac{\partial}{\partial x_i}|_p$, 
is the vector of components of $v$ with respect to coordinates $x$.
A Riemannian structure $g$ on M defines an inner product $<.,.>_p$ on $M_p$ such that 
$g_{ij}(x_1,...,x_n)=<\frac{\partial}{\partial x_i}|_p, \frac{\partial}{\partial x_j}|_p>_p$ are 
differentiable functions on $U$. At each $p\in {\bf x}(U)$, 
the $n\times n$ symmetric and positive definite matrix $G_x = \{g_{ij}(x)\}$ is called 
a coordinate representation of the metric at $p$.
If $y$ is another local system of coordinates at $p$ and 
$A = \{\frac{\partial y_j}{\partial x_i}|_p\}_{i=1,j=1}^{n,n}$, 
is the Jacobian of the change at $p$, which is a non-singular matrix, then 
component tangent vectors and metric representations change according to
\begin{equation}\label{covf_change_vector}
v_y = A v_x,
\end{equation}
\begin{equation}\label{covf_change_metric}
G_y = (A^{-1})' G_x A^{-1}.
\end{equation}

Any Riemmanian manifold can be endowed with a natural measure called {\it volume measure}.
Let $x$ be local coordinates at $p\in M$, the volume measure with respect to $x$ is defined by
$$
d\mathcal{V}(x) := (d\mathcal{V}(p))_x = \sqrt{det(G_x)} dx,
$$
where $dx$ is the Lebesgue measure in $\mathbb{R}^n$.
One easily checks that a change of local coordinates at p from $x$ to $y$, does not change 
the expression for $d\mathcal{V}(p)$, $d\mathcal{V}(x) = d\mathcal{V}(y)$.

Throughout this paper we will assume that M is a Riemannian manifold of dimension n.

{\subsection{Exponential map and its inverse}}

Geodesics on M are defined as solutions of first order system of differential equations, 
called geodesic equations (3.2 in \cite{carmo-riemannian}). 
In local chart $(U, {\bf x})$ they are
\begin{equation}\label{covf_geodesic_equations}
\frac{\partial x_k}{\partial t} = y_k, \\ 
\frac{\partial y_k}{\partial t} = - \sum_{i,j} \Gamma_{ij}^k y_iy_j, 
\end{equation}
where $\Gamma_{ij}^k$ are differentiable functions in $U$. The theory of 
ordinary differential equations says that for any $(x_1,y_1) \in U\times\mathbb{R}^n$, there 
exists a neighborhood $W$ of $(x_1,y_1)$ and $\epsilon >0$ such that 
for any $(x_0,y_0)\in W$, (\ref{covf_geodesic_equations}) has a unique solution 
$t\mapsto c(t)$ for $|t|<\epsilon$ satisfying $c(0)=x_0$ and $c'(0)=y_0$. Moreover, 
$c(t)$ depends differentially on the initial conditions.

For $q\in M$ and $v\in M_q$ let $\gamma(t, q, v)$, $t\in(-\epsilon,\epsilon)$ be 
a geodesic on M such that $\gamma(0, q,v)=q$ and $\gamma'(0, q,v)=v$.
Thus for $x(t)={\bf x}^{-1}\circ \gamma(t, q, v)$ and $y(t)=({\bf x}^{-1}\circ \gamma)'(t, q, v)$, $(x(t),y(t))$ is a solution of 
the system (\ref{covf_geodesic_equations}).

For any $p\in M$, there is a set $\mathcal{U}\subset TM$, $p \in \mathcal{U}$ and $\epsilon >0$, 
such that $\forall (q,v) \in \mathcal{U}$, 
$\gamma(t, q,v)$ is well defined and differentiable function of $(t,q,v)$ in $(-\epsilon,\epsilon+1) \times \mathcal{U}$. 
Let $\tilde q = \gamma(1,q,v)$ and $\tilde v = \gamma'(1,q,v)$. Then 
it follows from (\ref{covf_geodesic_equations}) that $(\tilde q, - \tilde v) \in \mathcal{U}$ and 
for $t\in(-\epsilon,\epsilon+1)$
\begin{equation}\label{covf_gamma_inverse}
\gamma(t,\tilde q, - \tilde v) = \gamma(1-t, q, v),
\end{equation}
$$
\gamma'(t,\tilde q, - \tilde v) = -\gamma'(1-t, q, v).
$$
The {\it exponential map}, $\exp: \mathcal{U} \to M$,  
is defined by
$$
\exp_q(v) = \exp(q,v) = \gamma(1,q,v).
$$
It is a differentiable map on $\mathcal{U}$.

For any $p\in M$, there is a maximal neighborhood $V(p)$ of the origin in $M_p$ 
where $\exp_p$ is a diffeomorphism; $U(p) = \exp_p(V(p))$ is called maximal normal neighborhood of $p$.
On $U(p)$, $\exp_p$ has an inverse, $$\exp_p^{-1}:U(p) \to V(p)\subset M_p,$$
which is also diffeomorphism.

The differential of the exponential map at $v\in V(p)$ 
$$
(d\exp_p)_v : M_p \to M_{\exp_pv}
$$
is an isomorphism.
One checks that $(d\exp_p)_O(v) = \gamma'(0, p, v) = v$ and 
\begin{equation}\label{covf_exp_v_v}
(d\exp_p)_v(v) = \gamma'(1, p, v).
\end{equation}
Indeed,
$$
(d\exp_p)_v(v) = \frac{d}{dt}\exp_p((t+1)v)|_{t=0} = 
\frac{d}{dt}\gamma(1,p,(t+1)v)|_{t=0} = \frac{d}{dt}\gamma(t+1,p,v)|_{t=0}.
$$
By the Gauss lemma (3.5 in \cite{carmo-riemannian}), $(d\exp_p)_v$ also satisfies 
$$
<(d\exp_p)_v (v), (d\exp_p)_v (w)> = <v,w>,
$$
for any $w\in V(p)$.

Let $p\in U(q)$ and $q\in U(p)$, then 
by applying (\ref{covf_exp_v_v}) and (\ref{covf_gamma_inverse}) we find 
$$
-\exp_q^{-1}p = -\gamma'(0,q,\exp_q^{-1}p) = 
$$
$$
\gamma'(1,p,\exp_p^{-1}q) = (d\exp_p)_{\exp_p^{-1} q} (\exp_p^{-1} q).
$$
Therefore, for a fixed $p$, we have following expression of $\exp_.^{-1} p : U(p) \to TM$ 
\begin{equation}\label{covf_exp_inverse}
\exp^{-1} p = - \gamma'(1,p, . ) \circ \exp_p^{-1}.
\end{equation}
Since $\gamma'(1,p, . )$ is differentiable in $V(p)$ and $\exp_p^{-1}$ is a diffeomorphism in $U(p)$, 
$\exp^{-1}p$ is differentiable in $U(p)$.

The map $q\mapsto \exp_q^{-1}p$ is differentiable in $U(p)$ in the following sense.
If $x$ are local coordinates at $q\in U(p)$, $q={\bf x}(x_1,...,x_n)$,  
then the components of $(\exp_q^{-1}p)_x \in \mathbb{R}^n$ are differentiable functions of $x$. 
Moreover, we show 
\newtheorem{covf_lemma1}{Lemma}
\begin{covf_lemma1}
For $q={\bf x}(x_1,...,x_n)\in U(p)$, the symmetric matrix 
\begin{equation}\label{covf_log_map_product}
Z_x(q, p) = (\exp_q^{-1}p)_x (\exp_q^{-1}p)_x'
\end{equation}
is differentiable in $x$.
\label{lemma:covf_lemma1}
\end{covf_lemma1}
A change of local coordinates from $x$ to $y$ with Jacobian $A$ at $q$, 
changes coordinate expression of $Z(q,p)$ according to 
\begin{equation}\label{covf_change_log_map_product}
Z_y(q, p) = A Z_x(q, p) A'.
\end{equation}

We will adopt, for brevity, following notation for the inverse exponential map 
$$
\overrightarrow{qp} := exp_q^{-1} p,
$$
in analogy to the Euclidean case where, 
$exp_q^{-1} p = p - q = \overrightarrow{qp}$, for $p,q \in \mathbb{R}^n$.

{\subsection{Tensors and tensor fields}}

Let $V$ be a n-dimensional vector space and $V^*$ be its {\it dual} space of linear functions on V.
Let also $x$ be a basis of V and $\tilde x$ be its dual basis.
For $v\in V$ with $v_x \in \mathbb{R}^n$ we denote 
the column vector of components of $v$, $v = \sum_i v_x^i x_i$.
While for a co-vector $w\in V^*$ with $w_x$ we denote the row vector 
of components of $w$, $w_x' \in \mathbb{R}^n$, $w = \sum_i w_x^i \tilde x_i$.

Co-variant 2-tensor $T$ on V is a bi-linear function 
$
T: V\times V \to \mathbb{R}, 
$
which with respect to the basis $x$ is represented by a matrix $T_x$. 
Coordinate expression for $T$ is 
$$
T(u,v) = u_x' T_x v_x, \\ \forall u,v \in V.
$$
With $T^2(V)$ we denote the vector space of co-variant 2-tensors on V.

Similarly, contra-variant 2-tensor $W$ is a bi-linear function
$
W: V^*\times V^* \to \mathbb{R}
$
with a coordinate expression $W_x$ with respect to $x$
$$
W(u,v) = u_x W_x v_x', \\ \forall u,v \in V^*.
$$
With $T_2(V)$ we denote the vector space of contra-variant 2-tensors on V.

Let $T\in T^2(V)$ and $W\in T_2(V)$. 
The contraction $TW$ of their tensor product $T\otimes W$ is a (1,1) tensor with
coordinates expression $T_xW_x$ with respect to $x$. 
We denote $TW\in T_1^1(V)$.

Let $y$ be another basis on $V$ such that $y =Ax$ for a non-singular matrix A.
Coordinate expressions for $T$, $W$ and $TW$ change according to 
$$
T_y = (A^{-1})'T_xA^{-1},
$$
$$
W_y = AW_xA',
$$
and
$$
(TW)_y = (A^{-1})'(TW)_xA'. 
$$
Recall that two matrices $C$ and $D$ are called congruent if 
there exists a non-singular matrix $P$ such that $C=PDP'$ 
and are called similar if $C=PDP^{-1}$.
Looking back at the change of coordinates rules we may conclude that the 
coordinate representations of (2,0) and (0,2) tensors are congruent, while 
those of (1,1) tensors are similar.

Let now M be a Riemmanian manifold with metric structure $G$.
The expression (\ref{covf_change_metric}) for the change of 
the metric representation at a point of M is co-variant like. 
Therefore, at each point $p\in M$ the metric is 
a symmetric and positive definite co-variant 2-tensor $G(p)$ on the tangent space $M_p$;
a fact that we notate with $G(p)\in T^2(M_p)$. Globally, $G$ is a co-variant 2-tensor field 
which is differentiable in the following sense. If $X,Y\in TM$ are 
two differentiable vector fields on M, then $G(X,Y)$ is a differentiable function on M.
With $\mathcal{T}^2(M)$ ($\mathcal{T}_2(M)$) we denote the differentiable co-variant (contra-variant) 
2-tensor fields on M. We write 
$$
G \in \mathcal{T}^2(M).
$$

Now we return back to $Z(q,p)$ which is given by
$
Z(q,p) = (\overrightarrow{qp})(\overrightarrow{qp})',
$
wherever the $\overrightarrow{qp}=\exp_q^{-1}p$ is defined. 
The change rule (\ref{covf_change_log_map_product}) for it is a contra-variant like 
and thus, $Z(q,p)$ is a symmetric and non-negative definite contra-variant 2-tensor at $M_q$. 
Moreover, by lemma \ref{lemma:covf_lemma1}, for any fixed $p\in M$, 
$Z(.,p)$ is a differentiable contra-variant tensor field on $U(p)$ - 
a fact that we write as 
\begin{equation}\label{covf_Z_tensor}
Z(.,p) \in \mathcal{T}_2(U(p)).
\end{equation}
Moreover, for every $q\in U(p)$, $G(q)Z(q,p) \in T_1^1(M_q)$ 
and $G(.)Z(.,p)$ is a differentiable (1,1)-tensor field on $U(p)$, i.e.
\begin{equation}\label{covf_GZ_tensor}
G(.)Z(.,p) \in \mathcal{T}_1^1(U(p)).
\end{equation}

{\subsection{Linear operators on tangent spaces}}

Linear operator on vector space $V$ is any $L:V \to V$ such that
$$
L(\alpha v_1 + \beta v_2) = \alpha L(v_1) + \beta L(v_2),
$$
for any two $v_1,v_2\in V$ and $\alpha,\beta\in \mathbb{R}$. 
With respect to a basis $x$, $L$ is represented by a matrix $L_x$ and then $L(v) = L_x v_x$.
Let $y$ be another basis such that $y=Ax$, then
$$
L(v) = L_y v_y = L_y Av_x = (L(v))_y = A(L(v))_x = AL_xv_x.
$$
Therefore
$$
L_y = AL_xA^{-1},
$$
which correspond to the change of coordinates rule for (1,1)-tensors.

Return back to Riemannian manifold setting. Let $p \in M$ and $q\in U(p)$, then 
$G(q)Z(q,p)$ defines a linear operator on $M_q$. In local coordinates $x$ at $q$ it is defined as 
$$
v_x\mapsto (G_xZ_x)' v_x  = Z_x G_x v_x\textrm{, } v\in M_q.
$$
If $v,w \in M_q$
$$
<w, (GZ)(v)> = w_x' G_xZ_xG_x v_x.
$$
In particular,
$$
<v, (GZ)(v)> = v_x' G_x Z_x G_x v_x = ((\exp_q^{-1}p)_x'G_x v_x)' ((\exp_q^{-1}p)_x'G_x v_x) > 0,
$$
for $v\ne 0$.

We summarize in the following 
\newtheorem{covf_lemma2}[covf_lemma1]{Lemma}
\begin{covf_lemma2}
For any fixed $p\in M$, $G(.)Z(.,p)$ is differentiable field of linear operators on $U(p)$ 
and thus, if $X$ is differentiable vector field on $U(p)$, 
then $(G(.)Z(.,p))(X)$ is also a differentiable vector field on $U(p)$.
\label{lemma:covf_lemma2}
\end{covf_lemma2}

{\subsection{Distributions on Riemannian manifolds}}

Let M be a Riemannian manifold and $(U_\alpha, {\bf x}_\alpha)$ are charts of M.
Open sets ${\bf x}_\alpha(U_\alpha)$ in M generate a $\sigma$-algebra on M which we will denote 
with $\mathcal{A}(M)$.
One easily verifies that the volume measure, $\mathcal{V}$, is a measure on the $\sigma$-algebra $\mathcal{A}(M)$ and thus,
$(M, \mathcal{A}(M), \mathcal{V})$ is a measure space.

A random variable X on M is any measurable function from a probability space 
$(\Omega,\mathcal{B},\mathcal{P})$ to $(M, \mathcal{A}, \mathcal{V})$.
The distribution function F of X is defined as 
$$
F(A) = \mathcal{P}(X^{-1}(A)),\\ A\in \mathcal{A}(M).
$$
F is a countably additive and satisfies $F\ge 0$, $F(\emptyset)=0$ and $F(M) = 1$.

\newtheorem{covf_def1}{Definition}
\begin{covf_def1}
A distribution F  on M is said to be absolute continuous with respect to 
the volume measure if 
$
F(A) = \int_A dF(p), \forall A\in \mathcal{A}(M),
$
where $dF(p)$ is given by
$$
dF(p) = f(p)dV(p), 
$$
for a $\mathcal{A}(M)$-measurable function $f$. 
We say that $f$ is density ({\it pdf}) of F.
\label{def:abscont_distr}
\end{covf_def1}
In the above definition, a density function $f$ is measurable in sense that $f^{-1}(B) \in \mathcal{A}(M)$  
for every {\it Borel} set B in $\mathbb{R}$. The density $f$ is continuous everywhere on M 
except eventually a set of volume measure zero.

In this work we will also consider discreet distributions on M. 
They are not absolute continuous w.r.t. $\mathcal{V}$ and instead of density have 
probability mass function ({\it pmf}).


{\section{Covariance fields}}

{\subsection{Definition}}

\newtheorem{covf_def2}[covf_def1]{Definition}
\begin{covf_def2}
Covariance field of probability distribution F is 
a contra-variant positive definite 2-tensor field $\Sigma$ on M, given by
$$
q \mapsto \Sigma(q) =  \int_{U(q)} (\overrightarrow{qp})(\overrightarrow{qp})' dF(p),
$$
where $U(q)$ is the maximal normal neighborhood of $q$.
\label{def:cov_field}
\end{covf_def2}
In the notation of (\ref{covf_log_map_product}), the covariance of F at $q$ is 
$$
\Sigma(q) =  \int_{U(q)} Z(q,p) dF(p) \in T_2(M_q)
$$
because $Z(q,p) \in T_2(M_q)$.
At this stage, 
we do not claim that $\Sigma$ is differentiable not even continuous field on M or on an open subset of M.

In local coordinates $x$, $Z_x(q,p)$ is a symmetric non-negative matrix and therefore 
$
\Sigma_x(q) =  \int_{U(q)} Z_x(q,p) dF(p)
$
is symmetric and non-negative definite. In fact, $\Sigma_x(q)$ will be positive definite, 
except the cases when the support of F in $x$ coordinates is a hyperplane in $\mathbb{R}^n$.
We ignore these cases, which obviously are caused by ill defined distributions, and assume 
positive definiteness of $\Sigma_x(q)$. Correspondingly, for the contra-variant tensor $\Sigma(q)$, 
we assume symmetry and positive definiteness.
The space of symmetric and positive definite matrices (tensors) we denote with $Sym_n^+$.

At point $q\in M$, we consider the product $G\Sigma(q) := G(q)\Sigma(q)$ of $G(q)\in T^2(M_q)$ and 
$\Sigma(q) \in T_2(M_q)$. It is a linear operator in $M_q$, i.e.
$$
G\Sigma(q) \in T_1^1 (M_q).
$$
Let $G_x$ and $\Sigma_x$ are representations of $G(q)$ and $\Sigma(q)$ 
with respect to coordinates $x$ about $q$. Then 
$$
<w, (G\Sigma(q))(v)> = w_x' G_x\Sigma_x G_x v_x,
$$
for $v,w\in M_q$. 
Moreover, $(G\Sigma(q))^{-1}$ is also a linear operator in $M_q$ and 
$$
<w, (G\Sigma(q))^{-1}(v)> = w_x' \Sigma_x^{-1} v_x,
$$

If $F(U(q)) = 1$, then  
$$
tr(G(q)\Sigma(q)) = \int_{U(q)} (\overrightarrow{qp})_x'G_x(q)(\overrightarrow{qp})_x dF(p) = 
$$
$$ 
\int_M ||Log_q p||^2 dF(p),
$$
Finally, for the intrinsic mean $\mu$, 
which is the Fr\'echet mean of F on the metric space M equipped with the geodesic distance, 
we obtain 
$$
\mu = argmin_q tr(G\Sigma(q)),
$$
in agreement to a well known fact about mean of distributions in $\mathbb{R}^n$.

After this illustrating example we are motivated to give a more general definition of expectation.
\newtheorem{covf_def3}[covf_def1]{Definition}
\begin{covf_def3}
Let $F$ be a distribution on M, $q\in M$, $TW \in T_1^1(M_q)$ and 
$h$ be a linear operator on $T_1^1(M_q)$.
Then the expectation of $h(TW)$ is defined to be 
$$
E(h(TW)) = \int_M h(TW(p)) dF(p).
$$
\label{def:expectation_tensor_invariant}
\end{covf_def3}
As shown above, this is applicable for $h(A) = tr(A)$ and $TW=G(q)Z(q,p)$. 
Then $G\Sigma(q) = E(tr(G(q)Z(q,p)))$, provided that $F(U(q))=1$.

{\subsection{Continuity of $G\Sigma$}}

We say that a series of points  $q_k$ on M converges to a point $q_0$ and denote $q_k \to q_0$ if 
$d(q_k,q_0) \to 0$, where $d$ is the geodesic distance on M.

\newtheorem{covf_prop1}{Proposition}
\begin{covf_prop1}\label{prop:covf_continuity_tr}
Let F be a distribution on M and 
for $q_k\in M$, $k=0,1,...$, we have $q_k \to q_0$, as $k \to \infty$, 
$F(B(q_k))=1$ for all $k$ and $tr(G\Sigma(q_0)) < \infty$.
Then 
$
tr(G\Sigma(q_k)) \to tr(G\Sigma(q_0)).
$
\end{covf_prop1}
{\it Proof.}
Without lost of generality we may assume that for all $q_k$, $d(q_0,q_k)\le r$. 
Let $X$ be a random variable with distribution F.
Define random variables $\xi_k = tr(G(q_k)Z(q_k,X))$ and $\xi_0 = tr(G(q_0)Z(q_0,X))$, where 
$Z(q,p)=(\overrightarrow{qp}) (\overrightarrow{qp})'$. Observe that 
$$
F(\cup_k(M\backslash B(q_k))) \le \sum_k (1-F(B(q_k))) = 0
$$ 
and therefore $F(\cap_k B(q_k)) = 1$.
Since $G$ and $Z$ are continuous at $q_0$ we have 
$
\xi_k \to \xi_0 \textrm{, a.e.}.
$
For every $p\in B(q_k)\cap B(q_0)$ 
$$
tr(G(q_k)Z(q_k,p)) = d^2(q_k,p) \le (d(q_k,q_0) + d(q_0,p))^2 
$$
$$
\le (r + d(q_0,p))^2 \le \max\{(r+1)^2, r^2 + (2r+1)d^2(q_0,p)\}
$$
and thus 
$$
0 \le \xi_k \le \max\{(r+1)^2, r^2 + (2r+1)\xi_0\} \textrm{, a.e.}.
$$
Finally, since $E\xi_0 < \infty$, the dominated convergence theorem gives us 
$
E(\xi_k) \to E(\xi_0),
$
which is exactly the claim.
\qquad$\Box$

\newtheorem{covf_prop2}[covf_prop1]{Proposition}
\begin{covf_prop2}\label{prop:covf_continuity_gsigma}
Under the conditions of proposition \ref{prop:covf_continuity_tr}, 
if $v_k\in M_{q_k} \to v_0\in M_{q_0}$ and $w_k\in M_{q_k} \to w_0\in M_{q_0}$, then 
$$
<v_k, (G\Sigma(q_k))(w_k)> \to <v_0, (G\Sigma(q_0))(w_0)>
$$
\end{covf_prop2}
{\it Proof.}
Define random variables 
$$\eta_k = <v_k, (G(q_k)Z(q_k,X))(v_k)>$$ and $$\eta_0 = <v_0, (G(q_0)Z(q_0,X))(v_0)>.$$
We have 
$
\eta_k \to \eta_0 \textrm{, a.e.}
$
and
$$
0 \le \eta_k \le ||v_k||tr(G(q_k)Z(q_k,X))||w_k|| .
$$
As in Proposition \ref{prop:covf_continuity_tr}, $\eta_k$ are bounded by a random variable with a finite expectation. 
Therefore, again by dominated convergence theorem, 
$
E(\eta_k) \to E(\eta_0).
$
\qquad$\Box$

We say that covariance field $\Sigma$ is continuous in the sense given by Proposition \ref{prop:covf_continuity_gsigma}, 
i.e. $G\Sigma$ is a continuous field of linear operators on tangent spaces of M.
\newtheorem{covf_def8}[covf_def1]{Definition}
\begin{covf_def8}
The covariance field $\Sigma$ of a probability distribution F is continuous at $q\in M$ if 
for any two continuous vector fields $v$ and $w$ on M defined in a neighborhood of $q$, the function 
$<v,(G\Sigma)w>$ is continuous at $q$.
\label{def:covf_cov_field_continuity}
\end{covf_def8}
Proposition \ref{prop:covf_continuity_tr} states sufficient conditions for continuity of $\Sigma$.
We emphasize that continuity of $G\Sigma$ may hold even when the density of F is discontinuous, 
provided the conditions are met. 
Moreover, even discrete distributions may have continuous covariance fields.

\newtheorem{covf_prop3}[covf_prop1]{Proposition}
\begin{covf_prop3}\label{prop:covf_continuity_gsigma_eigenvals}
Let F be a distribution on M, $q_0\in M$ such that $tr(G\Sigma(q_0)) < \infty$ 
and in a neighborhood $B_0$ of $q_0$ on M, $F(B(q)) = 1$, for all $q\in B_0$. 
Then the eignevalues of $G\Sigma$ are continuous functions at $q_0$.
\end{covf_prop3}
{\it Proof.}
Let $({\bf x}, U)$ be a parametrization about $q_0$, such that ${\bf x}(0) = q_0$ and 
${\bf x}(U)\subset B_0$. 
Let $\lambda_i(x)$, $i=1,...,n$, be the eigenvalues of $G\Sigma({\bf x}(x))$, $x\in U$.
Consider the canonical tangent vectors $v_i(x) = \frac{\partial}{\partial x_i}|_{{\bf x}(x)}$, which form 
continuous vector fields on ${\bf x}(U)$. By Proposition \ref{prop:covf_continuity_gsigma}, 
for all $i$ and $j$ 
$$
<v_i,(G\Sigma(q))(v_j)> = [G_x\Sigma_xG_x]_{ij} 
$$
are continuous functions at $0$. Since the elements of matrix $G_x$ are continuous at $0$, so 
are those of $G_x\Sigma_x$. Therefore the eigenvalues $\lambda_i(x)$ are continuous at $0$.
\qquad$\Box$


{\subsection{Extended covariance fields}}

In the course of our research we will find useful to extend the definition of covariance field 
to a whole class $\mathcal{COV}(F)$ of contra-variant tensor fields associated with a particular distribution $F$.
Its members are all 
\begin{equation}\label{def:covf_cov_extfield}
\Sigma(q; r) =  \int_M (\overrightarrow{qp})(\overrightarrow{qp})' r(||\overrightarrow{qp}||)dF(p),
\end{equation}
where $r:\mathbb{R}^+\to \mathbb{R}^+$ is a continuous function.
Thus,
$$
\mathcal{COV}(F) = \{ \Sigma(q; r)| r\in C(\mathbb{R}^+) \}.
$$
The role of $r$-function in (\ref{def:covf_cov_extfield}) is to control the 'amplitude', 
$tr(G\Sigma)$, of the covariance field. It has analytical purpose that is important in 
numerical experiments. Fields with large amplitude are difficult to be analysed numerically and 
this fact matters when one develops computational algorithms.
We will try to illustrate this with an example.
\newtheorem{covf_ex1}{Example}
\begin{covf_ex1}
Let $M = \mathbb{S}^2$ with the standard differential and metric structure 
(see details in Appendix S).
For every $p\in M$, $T_pM = \mathbb{R}^2$ and 
a normal neighborhood of $p$ is $S(p,\pi)$, the image of the circle $C(0,\pi)\subset T_pM$ 
under the exponential map, i.e. $S(p,\pi) = Exp_p(C(0,\pi))$.

Let F be uniform distribution on M and $y_1$,...,$y_k$ are samples drawn from F.
What is the sample covariance field of F based on these samples?

First, we consider the generic covariance field 
$$
\Sigma_1(y_j) = \sum_{i=1}^k (\overrightarrow{y_jy_i})(\overrightarrow{y_jy_i})',
$$
For a second one we apply $r(t) = (1-\frac{\pi}{2|t|})^2$, 
$$
\Sigma_2(y_j) = \sum_{i=1}^k (\overrightarrow{y_jy_i})(\overrightarrow{y_jy_i})'(1 - \frac{\pi}{2||\overrightarrow{y_jy_i}||})^2.
$$
Let $(\alpha,t)$ be polar coordinates on $\mathbb{S}^2$ at p. Then $dV(\alpha,t) = sin(t)d\alpha dt$ and the 
density is a constant, $f(\alpha,t) = 1/(4\pi)$.
As usual with $G$ we denote the metric tensor. First, we calculate
$$
E(tr(G(y_j)(\overrightarrow{y_jy_i})(\overrightarrow{y_jy_i})') = \frac{1}{2} \int_0^{\pi} t^2 sin(t) dt = \frac{\pi^2}{2} - 2,
$$
and
$$
E(tr(G(y_j)\Sigma_1(y_j))) = (k-1)(\frac{\pi^2}{2} - 2),
$$
while,
$$
E(tr(G(y_j)(\overrightarrow{y_jy_i})(\overrightarrow{y_jy_i})')(1 - \frac{\pi}{2||\overrightarrow{y_jy_i}||})^2 = 
$$
$$
\int_0^{\pi} t^2 sin(t) (1-\frac{\pi}{2t})^2dt = \frac{\pi^2}{4} - 2
$$
and
$$
E(tr(G(y_j)\Sigma_2(y_j))) = (k-1)(\frac{\pi^2}{4} - 2).
$$
Clearly, $E(tr(G(y_j)\Sigma_1(y_j))) > 6(k-1) E(tr(G(y_j)\Sigma_2(y_j)))$.
\newline
$\Box$
\end{covf_ex1}

Numerical experiments confirm the benefit of applying an amplitude bounding $r$-function.
A typical choice for $r$, suggested by experiments, is 
\begin{equation}
r(t) = (1-\frac{R}{2||t||})^2
\label{eq:cov_optimal_r}
\end{equation}
when the distribution has a bounded domain with geodesic radius $R$. 
It is a member of the family $\{r(t) = (1-\frac{a}{||t||})^2, a>0\}$.

Recall that a geodesic radius of distribution $F$ is the minimal $R$ such that 
for every $p\in supp(F)$, $supp(F) \subset Exp_p(C(R))$, 
where $C(R)$ is the ball with radius $R$ in tangent space at $p$, i.e. $C(R)=\{v\in\mathbb{R}^n| |v| \le R\}$.

Next result provides sharper estimator of parameter $a$ for the above family of $r$-functions 
based on the criterion $tr(G\Sigma)$ to be minimal.
But first we need following definition.

For every point $p\in M$, there are {\it geodesic spherical coordinates} $(\theta, t)$ defined by 
$$
p(\theta, t) = Exp_p(t\theta),\\ (\theta, t)\in U(p)\subset \mathbb{S}^{n-1}\times[0,\infty),
$$
and $Exp_p(U(p)) = \mathcal{B}(p)$, the maximal normal neighborhood of $p$.
The change of coordinates at $p$ from normal $v$ to spherical $(\theta, t)$ is 
$dv = t^{n-1}d\theta dt$.

A distribution $F$ is said to have a {\it bounded density} $f$, $dF(p) = f(p)dV(p)$, if 
there exists a positive constant $C$ such that 
for any $p\in M$, if $(\theta, t)$ are the spherical coordinates at $p$, then
$$
f(\theta, t)\sqrt{|G(\theta, t)|} \le C,\\ \forall \theta, t \in U(p).
$$
For example, on the unit 2-sphere, $\mathbb{S}^2$, 
$$
\sqrt{|G(\theta, t)|}d\theta dt = sin(t) d\theta dt, \\ (\theta,t)\in[0,2\pi)\times[0,\pi).
$$
and any function $f$ on $\mathbb{S}^2$ such that $C\ge f\ge 0$, after rescaling, defines a distribution with bounded density. 
\newtheorem{covf_lemma5}[covf_lemma1]{Lemma}
\begin{covf_lemma5}\label{covf_lemma:r_optimum_a}
Let $\mathcal{F}$ be a family of distributions on $M$ with geodesic radius $R$ and 
bounded by $C$ density, then for 
\begin{equation}
r(t; a) = (1-\frac{a}{||t||})^2,
\label{eq:cov_linear_r}
\end{equation}
we have 
$$
\sup_{F\in\mathcal{F}} \sup_{q\in M} tr(G(q)\Sigma(q;a)) \le 
\mathcal{V}(\mathbb{S}^{n-1}) C R (\frac{R^2}{3} - aR + a^2),
$$
where the covariance field $\Sigma(q,a)$ is defined by (\ref{def:covf_cov_extfield}) with $r(t,a)$ and 
$\mathcal{V}(\mathbb{S}^{n-1})$ is the volume of $\mathbb{S}^{n-1}$.
\label{lemma:cov_linear_r_estim}
\end{covf_lemma5}
{\it Proof.} For $q\in M$ and spherical coordinate system 
$(\theta, t)\in U\subset \mathbb{S}^{n-1}\times[0,R]$ at $p$ we have 
$$
tr(G(q)\Sigma(q;a)) = \int_M ||\overrightarrow{qp}||^2 (1 - \frac{a}{||\overrightarrow{qp}||})^2 dF(p) = 
$$
$$
\int_{U} t^2 (1 - \frac{a}{t})^2 f((\theta, t))\sqrt{|G((\theta, t))|} d\theta dt.
$$
From the bounded density assumption
$$
tr(G(q)\Sigma(q;a)) \le 
C \int_U t^2(1-\frac{a}{t})^2 d\theta dt \le 
$$
$$
\mathcal{V}(\mathbb{S}^{n-1}) C \int_{0}^{R} t^2(1-\frac{a}{t})^2 dt = 
\mathcal{V}(\mathbb{S}^{n-1}) C R (\frac{R^2}{3} - aR + a^2).
$$
Moreover the minimum of the right hand side function of $a$ is achieved for $R/2$.
$\Box$

As a consequence of lemma \ref{covf_lemma:r_optimum_a}, 
on $\mathbb{S}^2$, 
since all distributions have bounded geodesic radius $R=\pi$, 
$a=\pi/2$ is the optimal choice for the parameter in (\ref{eq:cov_linear_r}) and thus,  
$r(t) = (1-\frac{\pi}{2||t||})^2$ is the optimal member of the family (\ref{eq:cov_linear_r}).

{\section{Recovering distribution from covariance}}

In this section we consider the problem of recovering a distribution from its covariance field. 
If such a recovery is not possible, then the covariance structure will not be 
a complete distribution representation and 
its application potential will be diminished. 
Fortunately, in most circumstances, the answer of the problem is positive.

{\subsection{Similarity invariants}}

The space $Sym_n^+$ of symmetric and positive definite $n\times n$ matrices is a well studied manifold 
that accepts a Riemannian structure that makes it a symmetric space. 
For example, see A. Ohara, N. Suda, S. Amari (1996) and W. Forstner, B. Moonen (1999).
We define an important class of functions on $Sym_n^+$ that respects similarity operation and thus, 
are functions that can be applied on linear operators. 

\newtheorem{covf_def9}[covf_def1]{Definition}
\begin{covf_def9}\label{randcomp_def:siminvariant}
A {\it similarity invariant function} on $Sym_n^+$ is any continuous $h$ that satisfies
\begin{enumerate}
\item[(i)] $h(AXA',AYA') = h(X,Y)$, $\forall X,Y\in Sym_n^+$ and $A\in GL_n$.
\end{enumerate}
It is a non-negative with a unique root if
\begin{enumerate}
\item[(ii)] $h(X,Y) \ge 0$, $\forall X,Y\in Sym_n^+ \textrm{ and } h(X,Y) = 0 \iff X = Y$.
\end{enumerate}
Moreover, $h$ is called similarity invariant distance, if in addition to (i) and (ii) also satisfies
\begin{enumerate}
\item[(iii)] $h(X,Y) + h(Y,Z) \ge h(X,Z)$, $\forall X,Y,Z\in Sym_n^+$.
\end{enumerate}
We will denote with $\mathcal{SIM}_n$ the class of functions satisfying (i) and (ii).
\end{covf_def9}

Below we list several examples of similarity invariant functions. 

\begin{enumerate}
\item{For a fixed $Z\in Sym_n^+$, the similarity invariant 
$$
h_{trdif}(X, Y; Z) = |(tr(Z^{-1}X - Z^{-1}Y)|,
$$
satisfies (iii) but not (ii). Default choice will be $Z=G^{-1}$, the inverse of the metric tensor representation.
}
\item{The second one is sometimes referred as {\it affine-invariant distance} in $Sym_n^+$, 
see for example \cite{batchelor-tensors}, \cite{fletcher-tensors} and \cite{pennec-tensors}, 
and it is defined by 
$$
h_{trln2}(X, Y) = \{tr(ln^2(XY^{-1}))\}^{1/2}, X,Y\in Sym_2^+.
$$
Actually, $h_{trln2}$ is not a unique choice for a distance in $Sym_n^+$.
}
\item{
Log-likelihood function gives us another choice for h, 
$$
h_{lik}(X, Y) = tr(XY^{-1}) - ln|XY^{-1}| - n.
$$
It satisfies (i) and (ii) but it fails to satisfy the triangular inequality.
}
\item{
$$
h_{lntr}(X, Y) = ln( tr(XY^{-1} - YX^{-1})^2 ),
$$
$h_{lntr}$ is another candidate for a distance in $Sym_n^+$.
}
\item{
Another interesting choice for $h$ is
$$
h_{lnpr}(X, Y) = ln(tr(XY^{-1})tr(YX^{-1})),
$$
that satisfies (iii) and 'almost' satisfies (ii):  
$h_{lnpr}(X, Y) = 0 \iff X = cY$, for $c>0$.
}
\end{enumerate}

Let $F$ be a distribution on M and $\Sigma$ be its covariance field.
Consider a similarity invariant $h$ applied on covariance operator field $G\Sigma$. 
If the conditions of Proposition \ref{prop:covf_continuity_gsigma_eigenvals} hold, 
$h(G\Sigma)$ would be a continuous scalar field on M. Indeed, $h$ is 
a continuous function of the eigenvalues of $G\Sigma$, which are continuous by themselves. 

Scalar fields of the form $h(G\Sigma)$ can be viewed as {\it representations} of $F$. For some choices of $h$, 
they would be true distribution representations, in sense that the underlying distribution 
can be fully recovered from them. And this is what we are going to address next.

{\subsection{Recovering discrete distributions}}

Let $\mathcal{P}=\{p_i\}_{i=1}^{k}$ be a set of k points on M. 
By a discrete mass function ({\it pmf}) on M we understand any $f$ defined on the 
domain set $\mathcal{P}$, such that 
$f=\{f_i = f(p_i) \ge 0\}_{i=1}^k$ and $\sum_{i=1}^k f_i = 1$.
We write $f\in P_k^+$, where $P_k^+$ denotes the compact k-simplex. 

Let $\mathcal{Q}=\{q_j\}_{j=1}^{k}$ be another set of k points on M, called {\it observation set}, 
where the covariances of $f\in P_k^+$ will be considered. 

We assume that the set $\mathcal{P}$ is contained within the maximal normal neighborhood of 
each of the points $q_j$ in order for the vectors $\overrightarrow{q_jp_i}$ to be well defined. 
This assumption is not a strong one when M is a complete Riemannian manifold.
Thus, for every j we may assume a fixed local parametrization $(x_j,U_j)$ such that $\mathcal{P}\subset x_j(U_j)$.

Fix an amplitude controlling function $r$. Covariance of $f\in P_k^+$ at $q_j$ is defined as 
$$
\Sigma[f]_j := \Sigma[f](q_j) = 
\sum_{i=1}^k (\overrightarrow{q_jp_i})(\overrightarrow{q_jp_i})' r(||\overrightarrow{q_jp_i}||) f_i.
$$
Let us denote 
$$
Y_{ji} = (\overrightarrow{q_jp_i})(\overrightarrow{q_jp_i})' r(||\overrightarrow{q_jp_i}||),\\ i=1...k.
$$
then $\Sigma[f]_j = \sum_{i=1}^k f_i Y_{ji}$. The collection $\{\Sigma[f]_j\}_{j=1}^k$ is called a covariance set of $f$ on $\mathcal{Q}$. 

Now we are interested in the inverse problem. 
How to reconstruct a {\it pmf} from its observed covariances? 
Let $\mathcal{C}=\{C_j\in T_2(M_{q_j})\}_{j=1}^k$  be a set of contra-variant tensors, 
which may happen to be covariance tensors of an unknown distribution or may not. 
The problem is to find $f$ such that
$$
\Sigma[f]_j \approx C_j,\\ j=1...,k.
$$
To measure the 'closeness' we will use similarity invariant.

Define the functional 
\begin{equation}\label{covf_eq:H_function}
H(f) = \frac{1}{k} \sum_{j=1}^k h(\Sigma[f]_j, C_j), \\ f\in P_k^+, 
\end{equation}
where $h \in \mathcal{SIM}(n)$.
Now we can formulate more precisely our problem 
as an optimization one: find a {\it pmf} $\hat f$ such that
\begin{equation}\label{covf_eq:pmf_cov_estim1}
\hat f = argmin_f H(f).
\end{equation}
From the assumptions for $h$, it is clear that 
$$
H(f) = 0 \iff h(\Sigma[f]_j, C_j) = 0, \\ \forall j \iff \Sigma[f]_j = C_j, \\ \forall j.
$$

If $C_{j}$ are covariances that come from a {\it pmf} $f^0 \in P_k^+$, then $f^0$ will be a solution  of the system 
\begin{equation}\label{covf_eq:pmf_cov_equality}
\sum_{i=1}^k f_i Y_{ji} = C_j,\\ j=1...k.
\end{equation}
To be able to recover correctly $f^0$, the system (\ref{covf_eq:pmf_cov_equality}) 
should have a unique solution.
A necessary and sufficient condition for that is 
$$
rank(\mathcal{Y}|\mathcal{C}) = rank(\mathcal{Y}) = k, 
$$
where $\mathcal{Y} = \mathcal{Y}[q_j] := \{d_{ji}=tr(G(q_j)Y_{ji})\}_{j=1,i=1}^{k,k}$, $\mathcal{C} = \{tr(G(q_j)C_{j})\}_{j=1}^{k}$ 
and $\mathcal{Y}|\mathcal{C}$ is the matrix $\mathcal{Y}$ with vector $\mathcal{C}$ attached as a last column. 
Note that $d_{ji} = d^2(q_j,p_i)r(d(q_j,p_i))$ and for $r=1$ these are the squared geodesic distances.

\newtheorem{covf_def10}[covf_def1]{Definition}
\begin{covf_def10}\label{covf_def:covf_fullrank}
We say that a covariance operator field $G\Sigma$ on M has a full rank, if for any bounded subset $A\subset M$, 
$k\in\mathbb{N}$ and $k$-sample $\{q_j\}_{j=1}^k$, selected by a continuous distribution $Q$ on A, we have 
$P_{Q}(rank(\mathcal{Y}) = k) = 1$, where 
$\mathcal{Y} = \{d^2(q_j,p_i)r(d(q_j,p_i))\}_{i=1,j=1}^{k,k}$.
\end{covf_def10}

In Eulcidean space, $M\equiv\mathbb{R}^n$, the rank of $\{d^2(q_j,p_i)\}_{j=1,i=1}^{k,k}$ is bounded above by $n+2$, 
i.e. for a default covariance field in Euclidean space, 
$$
rank(\mathcal{Y}) \le n+2,
$$
and (\ref{covf_eq:pmf_cov_equality}) does not have a unique solution when $k>n+2$.
The problem can be fixed using a non-default covariance field with amplitude function $r\ne 1$ in (\ref{def:covf_cov_extfield}).

We hypothesize that on non-Euclidean space, a manifold with non-zero curvature, 
any covariance operator field is of full rank. Experiments on spheres, 
manifolds with constant sectional curvature +1, 
and hyperbolic plane, a manifold with constant sectional curvature -1, confirm 
the hypothesis, but of course a more formal argument is needed here.

If matrix $\mathcal{Y}$ has full rank, then one can find the {\it pmf} f directly solving system (\ref{covf_eq:pmf_cov_equality}) 
or minimizing the functional $H(f)$ as given in (\ref{covf_eq:H_function}) for $h\in\mathcal{SIM}_n$.
The second choice is much more general and gives us solutions even in cases 
when $C_j$ are not in fact true covariances, i.e. $rank(\mathcal{Y}|\mathcal{C}) > rank(\mathcal{Y})$.

It is also important to know in what cases the optimization problem (\ref{covf_eq:pmf_cov_estim1}) can be solved easily.
A function $h$ for which the corresponding functional $H$ is convex is an obvious choice since 
in that case it is straightforward to find the global minimum of $H$ by the gradient descend algorithm.

\newtheorem{covf_def11}[covf_def1]{Definition}
\begin{covf_def11}\label{def:recovdist_h_convex}
We say that $h\in\mathcal{SIM}_n$ is convex (in $Sym_n^+$) if for any $k$ and $Y_{ji},C_j\in Sym_n^+$, i,j=1,...,k, such that 
$rank(\mathcal{Y}) = k$, 
$\sum_{j=1}^k h(\sum_{i=1}^k f_iY_{ji}, C_j)$ is a convex function of $f$ in $P_k^+$.
\end{covf_def11}

From the list of invariants we list above, $h_{trdif}^2$, $h_{lik}$ and $h_{trsq}^2$ are convex.
We will show it for the last one.
\newtheorem{covf_ex2}[covf_ex1]{Example}
\begin{covf_ex2}
Let $Y_{ji}$ be such that $rank(\mathcal{Y} = \{tr(Y_{ji})\}_{j,i}) = k$. 
We will show that for $h_{trsq}^2(A,B) := tr((AB^{-1}-BA^{-1})^2)$, the functional 
$$
H(f) = \sum_{j=1}^k tr((\sum_{i=1}^k f_iY_{ij}) - (\sum_{i=1}^k f_iY_{ij})^{-1})^2
$$
is convex.
Without loss of generality we assume $C_i = I_n$.
Then 
$$
\frac{\partial H}{\partial f_s} = 2\sum_{j=1}^k tr[Y_{js}(\sum_{i=1}^k f_iY_{ij}) - (\sum_{i=1}^k f_iY_{ij})^{-1}Y_{sj}(\sum_{i=1}^k f_iY_{ij})^{-2})
$$
and defining $B_j = \sum_{i=1}^k f_iY_{ij}$, we obtain  
$$
\frac{\partial H}{\partial f_s} = 2\sum_{j=1}^k tr[Y_{js}(B_j-B_j^{-3})].
$$
Second derivatives are 
$$
\frac{\partial^2 H}{\partial f_s\partial f_l} = 2 \sum_j tr\{ Y_{sj}Y_{lj} + Y_{sj}B_j^{-1}Y_{lj}B_j^{-3} + Y_{sj}B_j^{-2}Y_{lj}B_j^{-2} + Y_{sj}B_j^{-3}Y_{lj}B_j^{-1} \}.
$$
For $w\in\mathbb{R}^k$ and $w\ne 0$, let $Z_j = \sum_{s=1}^k z_sY_{sj}$, then 
$$
w'\frac{\partial^2 H}{\partial f^2}w = 
2 \sum_{j=1}^k tr\{ Z_jZ_j + Z_{j}B_j^{-1}Z_{j}B_j^{-3} + Z_{j}B_j^{-2}Z_{j}B_j^{-2} + Z_{j}B_j^{-3}Z_{j}B_j^{-1} \}.
$$
Since $rank(\mathcal{Y})=k$, $\sum_{j=1}^k tr(Z_jZ_j) > 0$. Also, 
$$
tr(Z_{j}B_j^{-1}Z_{j}B_j^{-3}) = tr([B_j^{-2}Z_{j}B_j^{-1}][B_j^{-1}Z_{j}B_j^{-2}]') \ge 0
$$
and similarly $tr(Z_{j}B_j^{-2}Z_{j}B_j^{-2})\ge 0$.
That proves the convexity of $h_{trsq}^2$.
$\Box$
\end{covf_ex2}

Now we are interested in the problem of consistency of estimators (\ref{covf_eq:pmf_cov_estim1}).
If we assume that the covariances $C_j$ are random and converge in probability to some matrices, 
is it true that the estimators $\hat f$ also converge? To guarantee a positive answer we need more assumptions for the  
invariant $h$. 

\newtheorem{recovdist_th1}{Theorem}
\begin{recovdist_th1}\label{th:recovdist_discrete}
Let $h\in\mathcal{SIM}_n$ be a distance function. 
Let also $f^0\in P_k^+$ be a {\it pmf} on $\mathcal{P}$ and $\{C_j^0 = \Sigma[f^0](q_j)\}_{j=1}^k$ 
be its covariance set on $\mathcal{Q}$.
If $\mathcal{C}^m=\{C_j^m\in Sym_n^+\}_{j=1}^k$ is a sequence of set of random matrices such that 
$$
\forall j = 1,...,k, ~~ h(C_j^m,C_j^0) \longrightarrow_p 0,\textrm{ as }m\to\infty
$$
then for 
$$
\hat f^m = argmin_{f\in P_k^+} H_m(f),\textrm{ where } H_m(f) = \frac{1}{k} \sum_{j=1}^k h(\Sigma[f]_j, C_j^m),
$$
we have $\hat f^m \longrightarrow_p f^0$.
\end{recovdist_th1}
{\it Proof.}\\
Since $h$ satisfies the triangular inequality, then for any $f\in P_k^+$
$$
|h(\Sigma[f]_j, C_j^m) - h(\Sigma[f]_j, C_j^0)| \le h(C_j^m,C_j^0) 
$$
and
$$
sup_{f\in P_k^+}|h(\Sigma[f]_j, C_j^m) - h(\Sigma[f]_j, C_j^0| \le h(C_j^m,C_j^0).
$$
By summing on j we obtain
$$
sup_{f\in P_k^+}|H_m(f) - H_0(f)| \le \sum_{j=1}^k h(C_j^m,C_j^0).
$$

Under the assumptions on $h$, $H_0(f)$ has a well-separated minimum  
at $f^0\in P_k^+$. In fact $H_0(f^0) = 0$.
Therefore, for any $\delta > 0$, there exists $\epsilon > 0$, such that
$$
H_0(f) > H_0(f^0) + 2\epsilon, \textrm{ for any f such that } |f - f^0|_{L_2} > \delta.
$$
Then we have
$$
P(|\hat f^m - f^0|_{L_2} > \delta) \le P(H_0(\hat f^m) > H_0(f^0) + 2\epsilon) \le 
$$
$$
P(H_0(\hat f^m) - H_m(\hat f^m) + H_m(\hat f^m) - H_0(f^0) > 2\epsilon ) \le 
$$
since $H_m(\hat f^m) \le H_m(f^0)$
$$
P(H_0(\hat f^m) - H_m(\hat f^m) + H_m(f^0) - H_0(f^0) > 2\epsilon ) \le 
$$
$$
P(2 sup_{f\in P_k^+}|H_m(f) - H_0(f)| > 2\epsilon) \le
$$
$$
P( \sum_{j=1}^k h(C_j^m,C_j^0) > k\epsilon) \le \sum_{j=1}^k P( h(C_j^m,C_j^0) > \epsilon).
$$
Since for any j, $P( h(C_j^m,C_j^0) > \epsilon) \longrightarrow 0$, we have 
$P(|\hat f^m - f^0|_{L_2} > \delta) \longrightarrow 0$. $\Box$

Unfortunatelly, the above result is not very useful in practice for distance functions are usually non-convex and 
non-convexity of $H$ makes its optimization difficult. 
In fact, the condition on $h$ to satisfy the triangular inequality is a stronger assumption than what we actually need. 
We observe that it is only used to bound uniformly 
$|h(\Sigma[f]_j, C_j^m) - h(\Sigma[f]_j, C_j^0|$ by 
$h(C_j^m, C_j^0)$. Therefore, if we guarantee the uniform convergence of the former difference by other means, 
the triangular inequality condition will be redundant. 

\newtheorem{recovdist_cons1}{Consistency Criterion}
\begin{recovdist_cons1}
We say that a similarity invariant function $h$ satisfies the consistency criterion, if 
for $C^m\in Sym_n^+$ random, and $B_i,C^0\in Sym_n^+$, 
$$
h(C^m,C^0) \longrightarrow_p 0 
$$
implies
$$
\sup_{f} |h(\sum_{i=1}^k f_iB_i, C^m) - h(\sum_{i=1}^k f_iB_i, C^0)| \longrightarrow_p 0, \textrm{ as } m\to\infty.
$$
\label{def:recovdist_consistency_criterion1}
\end{recovdist_cons1}
The following theorem is a stronger version of Theorem \ref{th:recovdist_discrete}, but using the above consistency criterion, 
and can be proven similarly.
\newtheorem{recovdist_th2}[recovdist_th1]{Theorem}
\begin{recovdist_th2}\label{th2:recovdist_discrete}
Let $h\in\mathcal{SIM}_n$ satisfy consistency criterion \ref{def:recovdist_consistency_criterion1}.
Let also $f^0\in P_k^+$ be a {\it pmf} on $\mathcal{P}$ and $\{C_j^0 = \Sigma[f^0](q_j)\}_{j=1}^k$ 
be its covariance set on $\mathcal{Q}$.
If $\mathcal{C}^m=\{C_j^m\in Sym_n^+\}_{j=1}^k$ is a sequence of set of random matrices such that 
$$
\forall j = 1,...,k, ~~ h(C_j^m,C_j^0) \longrightarrow_p 0,\textrm{ as }m\to\infty
$$
then $\hat f^m \longrightarrow_p f^0$.
\end{recovdist_th2}

It turns out that $h_{trln2}$ invariant is a distance but is not convex in the sense of 
definition (\ref{def:recovdist_h_convex}) and finding the global minimum of $H_{trln2}$ is difficult. 
On the other hand invariants $h_{lik}$ and $h_{trsq}^2$ are both convex and 
satisfy the consistency criterion \ref{def:recovdist_consistency_criterion1}, which makes them better choices.

The condition $h(C_j^m,C_j^0) \longrightarrow_p 0$ can be further simplified if $h$ is continuous. 
We say that a sequence $X_m$ of random $n\times n$ matrices converges in probability to matrix $C$ and write 
$X_m\longrightarrow_p C$,if for any $v\in \mathbb{R}^n$, $v'(X_m-C)v \to_p 0$.
One can easily check that if $h(X,C)$ is continuous in $X\in Sym_n^+$ for every $C\in Sym_n^+$, then 
$X_m\longrightarrow_p C$ if and only if $h(X_m, C) \to_p 0$. 
Next we have the following corollary of Theorem \ref{th2:recovdist_discrete}.
\newtheorem{recovdist_cor1}{Corollary}
\begin{recovdist_cor1}
Let $h\in\mathcal{SIM}_n$ be a continuous invariant that satisfies consistency criterion \ref{def:recovdist_consistency_criterion1}.
If 
$C_j^m \longrightarrow_p C_j^0$, $\forall j = 1,...,k$, 
then $\hat f^m \longrightarrow_p f^0$.
\end{recovdist_cor1}
Note that both $h_{lik}$ and $h_{trsq}^2$ satisfy the conditions of the above corollary.

{\subsection{Recovering continuous distributions}}

We will give a constructive procedure for recovering a continuous distribution density 
from its covariance operator field. 
It turns out that any covariance field of full rank specifies completely the underlying density 
when its domain is a bounded compact. This fact shows that covariance fields are in general 
faithful representations of corresponding distributions. Although this fact may seem obvious at first look, 
there is a notable exception that makes the problem of recovering relevant.
In Euclidean space, $M=\mathbb{R}^n$, we have the following relation for the default covariance field 
of random variable $X$ with $E(X) = \mu$
$$
E[(q-x)(q-x)'] = E[(p-x)(p-x)'] + (q-p)(q-p)' + (q-p)(p-\mu)' + (p-\mu)(q-p)',
$$
Thus 
$$
\Sigma(q) = \Sigma(\mu) +  (q-\mu)(q-\mu)'
$$
does not contain any information beyond the first two moments: the mean $\mu$ and the covariance $\Sigma(\mu)$.
Therefore the covariance field is defined only by the first two moments of $X$ and 
can not possibly represent the whole distribution.

We use the same approach as for recovering discrete distributions and start with 
selecting appropriate similarity invariants as technical instruments. 
We need to make stronger assumptions than that in consistency criterion (\ref{def:recovdist_consistency_criterion1}).
\newtheorem{recovdist_cons2}[recovdist_cons1]{Consistency Criterion}
\begin{recovdist_cons2}
We say that similarity invariant function $h$ satisfies the consistency criterion if 
for $B_i, C_1, C_2\in Sym_n^+$ such that $||B_i||\le \gamma$ and $||C_i||\le \gamma$
$$
\sup_{f} |h(\sum_{i=1}^k f_iB_i, C_1) - h(\sum_{i=1}^k f_iB_i, C_2)| \le 
\alpha\gamma ||C_1 - C_2||,
$$
for a constant $\alpha>0$, independent of $B_i$ and $C_i$.
\label{def:recovdist_consistency_criterion2}
\end{recovdist_cons2}
\newtheorem{covf_ex3}[covf_ex1]{Example}
\begin{covf_ex3}
We will show that $h_{trdif}^2(A,B) = tr^2(A-B)$ is convex and satisfies the consistency criterion \ref{def:recovdist_consistency_criterion2}.
Let $Y_{ji}$ are such that $rank(\mathcal{Y} = \{tr(Y_{ji})\}_{j,i}) = k$.
Define $L(f) = \sum_{j=1}^k (\sum_{i=1}^k f_i tr(Y_{ij}) - tr(C_j))^2$, then 
$$
\frac{\partial L}{\partial f_s} = 2\sum_{j=1}^k tr(Y_{js})(\sum_{i=1}^k f_i tr(Y_{ij}) - tr(C_j))
$$
and 
$$
\frac{\partial^2 L}{\partial f_s\partial f_l} = 2\sum_{j=1}^k tr(Y_{js})tr(Y_{jl})
$$
For $w\in\mathbb{R}^k$, 
$$
w'\frac{\partial^2 L}{\partial f^2}w = 2 \sum_{j=1}^k (\sum_{s=1}^k w_s tr(Y_{js}))^2 > 0
$$
because $rank(\mathcal{Y})=k$. This shows the convexity. 

Next, observe that for $B_i$, $C_1$ and $C_2$, such that $||B_i||\le \gamma$, $||C_i||\le \gamma$
$$
|(\sum_{i=1}^k f_i tr(B_i) - tr(C_1))^2 - (\sum_{i=1}^k f_i tr(B_i) - tr(C_2))^2| \le 
$$
$$
|tr(C_1) - tr(C_2)|~|2tr(\sum_{i=1}^k f_i B_i) + tr(C_1) + tr(C_2)| \le 2n(n+1)\gamma||C_1-C_2||,
$$
since for any $n\times n$ matrix $X$, $|tr(X)|\le n||X||$.
$\Box$
\end{covf_ex3}

In order to show our main result we need to put some restrictions on the size of distribution domains. 
In the sequel, $d_g$ denotes the geodesic distance corresponding to a metric $g$ on M. 
\newtheorem{covf_def12}[covf_def1]{Definition}
\begin{covf_def12}\label{def:recovdist_lipschitz}
We say that the inverse exponential map $exp^{-1}$ on M is Lipschitz in $A\subset M$, 
if there exists a constant $\beta > 0$ such that $\forall q\in A$ 
$$
||\overrightarrow{qp_1} - \overrightarrow{qp_2} ||_{L^2} = ||exp_q^{-1}p_1 - exp_q^{-1}p_2||_{L^2} \le \beta d_g(p_1,p_2),\\ \forall p_1,p_2\in A.
$$
\end{covf_def12}
Necessary condition for $exp^{-1}$ to be Lipschitz in $A$ is for any $q\in A$, $A\subset U(q)$, 
where, recall, $U(q)$ denotes the domain where $exp_q^{-1}$ is defined and is in fact diffeomorphic. 
It is an open question whether this condition is sufficient one. The next example provides some evidence in its support.

\newtheorem{covf_ex4}[covf_ex1]{Example}
\begin{covf_ex4}
On the unit 2-sphere, $M=\mathbb{S}^2$, the $exp^{-1}$-map is Lipschitz in any compact $K$ with 
geodesic diameter less than $\pi$. If $\rho=diam(K) < \pi$ then parameter $\beta = \frac{\rho}{sin(\rho)}$.
\end{covf_ex4}

\newtheorem{recovdist_th3}[recovdist_th1]{Theorem}
\begin{recovdist_th3}\label{th:recovdist_dcontinuous}
Let $F$ be a distribution on M and let $K$ be a compact in M, which is bounded, $diam(K)\le R$, 
contains the support of $F$, $F(K) = 1$, 
and $Log$-map is Lipschitz in $K$.
Let also the default covariance field $\Sigma$ of $F$ has a full rank $a.e.$ on M.
Then there is a sequence $\hat F^m$ of {\it pmfs} on M obtained from the field $G\Sigma$ alone, 
such that $\forall V \subset M$
$$
\hat F^m(V) \longrightarrow F(V),\textrm{ as }m\to\infty.
$$
\end{recovdist_th3}
{\it Proof.}
Under the assumptions for $K$, for any $m > 0$ there exists a partition 
$\{U_j^m\}_{j=1}^{N(m)}$ of $K$ with geodesic diameter no greater than $1/m$ and 
each point $q_j^m \in U_j^m$ is selected independently by the uniform distribution on $U_j^m$.
Then $P(rank(\mathcal{Y}[q_j^m]) = N(m)) = 1$. 

Define $\tilde f_j^m = F(U_j^m)$ and 
$$
\widetilde \Sigma^m(q) = \sum_{j=1}^{N(m)} \tilde f_j^m (\overrightarrow{q q_j^m})(\overrightarrow{q q_j^m})',\\ q\in K.
$$
Note that $\widetilde \Sigma^m(q)$ is random because it depends on the choice of $q_j^m$.
Since 
$$
\sum_{j=1}^{N(m)}  (\overrightarrow{q q_j^m})(\overrightarrow{q q_j^m})' - 
\int_{K} (\overrightarrow{q p})(\overrightarrow{q p})' dF(p) = 
$$
$$
\sum_{j=1}^{N(m)} \int_{U_j^m}[(\overrightarrow{q q_j^m})(\overrightarrow{q q_j^m})' - (\overrightarrow{q p})(\overrightarrow{q p})'] dF(p), 
$$
for any $v\in M_q$ we have
$$
|<v, (G\widetilde\Sigma^m(q) - G\Sigma(q))(v)>| = 
$$
$$
\sum_{j=1}^{N(m)} \int_{U_j^m} |<v,\overrightarrow{q q_j^m}>^2 - <v,\overrightarrow{q p}>^2| dF(p) \le 
$$
$$
\sum_{j=1}^{N(m)} \int_{U_j^m} ||v||~||\overrightarrow{q q_j^m} - \overrightarrow{q p}||~||v||~(||\overrightarrow{q q_j^m}|| + ||\overrightarrow{q p}||) dF(p) < 
$$
$$
||v||^2 \beta (2R^2 + \frac{1}{m})\frac{1}{m}, \textrm{ w. p. } 1.
$$
Thus the operator norm of $G\widetilde\Sigma^m(q) - G\Sigma(q)$ is uniformly bounded
\begin{equation}\label{recovdist_eq:tildesigma_uniform_bound}
||G\widetilde\Sigma^m(q) - G\Sigma(q)|| \le \beta (2R^2 + \frac{1}{m})\frac{1}{m}, \textrm{ w. p. } 1.
\end{equation}
In particular, since dim(M)=n 
$$
\max_{q\in A} |tr(G\widetilde\Sigma^m(q) - G\Sigma(q))| \le  n\beta (2R^2 + \frac{1}{m})\frac{1}{m}, \textrm{ w. p. } 1.
$$

Since $\Sigma$ has a full rank $a.e.$ on M, then there exists a $h\in\mathcal{SIM}_n$ 
that is convex (in sense of definition \ref{def:recovdist_h_convex}) and 
satisfies consistency criteria \ref{def:recovdist_consistency_criterion2} (w.p. 1). 
For example, we may take $h=h_{trdif}^2$. 
Define
$$
H_m(f) = \frac{1}{N(m)} \sum_{j=1}^{N(m)} h(G\Sigma[f](q_j^m), G\Sigma(q_j^m))
$$
$$
\tilde H_m(f) = \frac{1}{N(m)} \sum_{j=1}^{N(m)} h(G\Sigma[f](q_j^m), G\tilde\Sigma^m(q_j^m))
$$
and let 
$$
\hat f^m = argmin_{f} H_m(f) \textrm{ and } \tilde f^m = argmin_{f} \tilde H_m(f).
$$

Since $\tilde H_m(\tilde f^m)=0$, $\tilde H_m(f)$ has a well separated minimum at $\tilde f^m$, i.e. 
for any $\delta > 0$, there exists $\epsilon > 0$, such that
$$
\tilde H_m(f) > \tilde H_m(\tilde f^m) + 2\epsilon = 2\epsilon, \textrm{ for any }f\textrm{ such that } |f - \tilde f^m|_{L_2} > \delta.
$$
Therefore 
$$
P(|\hat f^m - \tilde f^m|_{L_2} > \delta) \le P(\tilde H_m(\hat f^m) > \tilde H_m(\tilde f^m) + 2\epsilon) = 
$$
$$
P(\tilde H_m(\hat f^m) - H_m(\hat f^m) + H_m(\hat f^m) - \tilde H_m(\tilde f^m) > 2\epsilon ) \le 
$$
since $H_m(\hat f^m) \le H_m(\tilde f^m)$
$$
P(\tilde H_m(\hat f^m) - H_m(\hat f^m) + H_m(\tilde f^m) - \tilde H_m(\tilde f^m) > 2\epsilon ) \le 
$$
$$
P(2 \sup_{f\in P_k^+}|H_m(f) - \tilde H_m(f)| > 2\epsilon) \le
$$
by consistency criteria \ref{def:recovdist_consistency_criterion2} for $h$ and 
the fact that $||(\overrightarrow{q p})(\overrightarrow{q p})'||\le nR^2$ on $K$
$$
P( \alpha nR^2 \frac{1}{N(m)} \sum_{j=1}^{N(m)} ||G\widetilde\Sigma^m(q_j^m) - G\Sigma(q_j^m)|| > \epsilon) \longrightarrow 0,\textrm{ as }m\to\infty 
$$
by the virtue of (\ref{recovdist_eq:tildesigma_uniform_bound}).
Thus
$
\hat f^m - \tilde f^m \longrightarrow_p 0.
$

Finally, for any $V\subset M$ we have 
$$
\sum_{j:q_j^m\in V} \tilde f_j^m \longrightarrow_p F(V), \textrm{ as } m\to\infty
$$
and therefore 
$$
\hat F^m(V) = \sum_{j:q_j^m\in V} \hat f_j^m \longrightarrow_p F(V).
$$
$\Box$

The above theorem is constructive and give us the freedom to choose 
a similarity invariant $h$ provided it is convex and satisfies the consistency criterion \ref{def:recovdist_consistency_criterion2}. 
As we showed, $h_{trdif}^2$ satisfies both requirements and can be applied.
In this case the corresponding optimization functional
$$
H_m(f) = \frac{1}{N(m)} \sum_{j=1}^{N(m)} (\sum_i d^2(q_i^m,q_j^m) - \rho(q_j^m))^2
$$
is based on the function $\rho(q) = E_F d^2(q,p) = \int_K d^2(q,p)dF(p)$, the "variance" of $F$ with respect to $q$. 
What the theorem states, basically, is that distribution $F$ can be recovered from the scalar field $\rho$ on M 
provided that the full-rank condition for the covariance field on M holds and in the domain of $F$, $exp^{-1}$-map is Lipschitz. 

The requirement for $h$ to satisfy the consistency criterion \ref{def:recovdist_consistency_criterion2} can be relaxed. 
We only need $h\in\mathcal{SIM}_n$ to be continuous invariant that satisfies consistency criterion \ref{def:recovdist_consistency_criterion1}. 
Indeed, looking back at the proof above we see that 
$||G\widetilde\Sigma^m(q_j^m) - G\Sigma(q_j^m)|| \longrightarrow_p 0$, and therefore 
$\sup_{f\in P_k^+}|H_m(f) - \tilde H_m(f)| \longrightarrow_p 0$, which is what we needed to show 
$\hat f^m - \tilde f^m \longrightarrow_p 0$. 
That observation gives us more choices for $h$ such as the convex invariants $h_{lik}$ and $h_{trsq}^2$.

In $\mathbb{R}^n$ recovering from default covariance field is not possible because the full-rank condition fails. 
On non-Euclidean (with non-zero curvature) spaces however, the full-rank condition is generally true and reconstruction is possible. 
For example, if $M=\mathbb{S}^2$ and $supp(F) \subset K$, for a compact on $\mathbb{S}^2$ with $diam(K) < \pi$, then 
the theorem is applicable and $F$ can surely be recovered.

Some of the conditions of the above theorem can be relaxed. For example, 
as defined the full-rank condition for the covariance field allow infinitely many choices 
for the points $q_j^m$ and thus infinitely many choices of sequence $\hat F^m$ converging to $F$. 
Different invariants $h$ also give different approximating distributions $\hat F^m$.

Theorem \ref{th:recovdist_discrete} shows how to recover discrete distribution 
provided the full-rank condition only.
Lipschitz condition in theorem \ref{th:recovdist_dcontinuous} is a technical one and it might be possible to relax it in 
a different approach to the problem.
Another opportunity, for example, is to work with a covariance field with amplitude $r\ne 1$ and 
to satisfy both full-rank condition and Lipschitz condition, even in $\mathbb{R}^n$.

{\subsection{Summary}}

We introduced the concept of covariance field of a distribution on a Riemannian manifold. 
It is a contra-variant 2-tensor field on the manifold. 
Closely associated with it is a covariance operator field, which defines a linear operator on the tangent space 
at each point on the manifold. Covariance operator fields, in most cases, are continuous. 

Covariance fields, in general, recover the underlying distributions and in this sense are faithful distribution representations. 
There is one major exception though where such reconstruction is not possible, the Euclidean space $\mathbb{R}^n$. 
This interesting fact shows that the covariance field concept is indeed more relevant to non-Euclidean manifolds. 


\end{document}